\documentclass[11pt]{amsart}
\usepackage{fullpage}

\def\C{{\hbox{\bf C}}}

\def\P{{\hbox{\bf P}}}
\def\E{{\hbox{\bf E}}}
\def\V{{\hbox{\bf Var}}}

\def\det{{\hbox{det}}}
 at 10 true pt

\def\eps{\varepsilon}

\def\per{{\hbox{\rm per}}}

\def\spec{\hbox{\rm spec} }

\def \det {{\rm {det}}}

\theoremstyle{plain}
\newtheorem{theorem}[subsection]{Theorem}
\newtheorem{conjecture}[subsection]{Conjecture}

\newtheorem{lemma}[subsection]{Lemma}
\newtheorem{corollary}[subsection]{Corollary}

\newtheorem{question}[subsection]{Question}
\theoremstyle{remark}
\newtheorem{remark}[subsection]{Remark}

\theoremstyle{definition}
\newtheorem{definition}[subsection]{Definition}

\begin{document}
\title{Concentration of random determinants and permanent estimators}

\author{Kevin P. Costello}
\address{School of Mathematics, Georgia Institute of Technology, Atlanta GA 30308}
\email{kcostell@math.gatech.edu}
\thanks{K. Costello was supported in this research by NSF Grant DMS-0635607}
\author{Van Vu}
\address{Department of Mathematics, Rutgers University, Piscataway NJ 08854}
\email{vanvu@math.rutgers.edu}
\thanks{V. Vu  is supported by  NSF Career Grant 0635606.}
\begin{abstract}

We show that the absolute value of the determinant of a matrix with random
 independent (but not necessarily iid) entries is strongly concentrated around its mean.

As an application, we show that
 Godsil-Gutman and Barvinok estimators for the permanent of a strictly positive matrix give sub-exponential approximation ratios with high probability.

 A positive answer to the main conjecture of the paper would lead to polynomial approximation ratios in the above problem.
\end{abstract}
\maketitle

\section{Introduction}
Let $A$ be an $n \times n$ square  matrix. We denote by $\det A$ and $\per A$
its determinant and permanent, respectively, which are defined by
\begin{equation*}
 \det A=\sum_{\sigma} (-1)^{sgn \sigma} \prod_{i=1}^n a_{i\sigma i}, \, \, \, \, \, \, \, \, \, \, \per A=\sum_{\sigma} \prod_{i=1}^n a_{i\sigma i},
\end{equation*}
where the sum is taken over all permutations in $S_n$ and $a_{ij}$ denotes the $(i,j)$ entry of $A$.  

 In this paper, we focus on a random matrix $A$  whose entries are independent (but not necessarily iid) random variables with mean zero. The size of $A$, (which we denote by $n$) should be thought of as tending to infinity and all asymptotic notation will be used under this assumption.

Our main concern is the following basic question

\begin{question} \label{question:main}
How is  $|\det A|$  distributed ?
\end{question}
A special  case is when the entries of $A$ are iid Gaussian (with variance one). In this case,
it is known  that $\log |\det A|$ satisfies the central limit theorem.

\begin{theorem} \label{theorem:normal}
Let $A$ be the random matrix of size $n$ whose entries are iid Gaussian with variance one. Then
\begin{equation*}
 \frac{\log(|\det A|)- \frac{1}{2}\log( (n-1)!)}{\sqrt{\frac{\log n}{2}}}
\end{equation*}
converges weakly to the standard Gaussian variable $N(0,1)$.
\end{theorem}

This statement is easy to verify, as one can write

$$|\det A| =\prod_{i=1}^{n} d_i $$

where $d_i$ is the distance from the $i$th row vector of $A$ to the subspace spanned by the first $i-1$ rows.
As $A$ has iid Gaussian entries, the random variables $d_i$ are independent. Furthermore, their distributions
can be computed explicitly and the theorem follows from  Lyapunov's Central Limit Theorem
and a routine calculation. (We include the details in Appendix A for the reader's convenience.)

The situation with  general random matrices is considerably more complicated. In \cite{Gir},
Girko claimed that Theorem \ref{theorem:normal} still holds if
the entries are no longer Gaussian, but still iid with mean zero and variance one.
We believe that this statement is true, but could not understand Girko's proof.
 On the other hand, it seems possible  that one can give an alternative proof using recent developments in the field.

In this paper, instead of limiting distribution, we focus on tail inequalities, which are usually very
 useful in probabilistic combinatorics and related fields. As an illustration, we present an  application concerning the problem of
computing the permanent using determinant estimators. A consequence of our main result shows that  one can
use a determinant estimator to estimate the permanent of a  matrix of size $n$ with positive entries within a sub-exponential factor $\exp(n^{2/3})$ with high probability. If Conjecture \ref{conj:conj1} holds, then the approximation will typically be within a polynomial factor $n^{O(1)}$.

To start, we note an old observation of Tur\'an  that if the entries
of $A$ are iid with mean 0 and variance 1, then $\E(|\det A|^2) =
\E(\det A^2)=n!$. Combining this with Theorem \ref{theorem:normal},
we obtain the following corollary.

\begin{corollary} \label{cor:gaussian}
Let $A$ be the random matrix of size $n$ whose entries are iid Gaussian with variance one. Then
with probability tending to one

\begin{equation}
|\det A^2| = n^{-1+o(1)} \E (\det A^2).
\end{equation}
\end{corollary}

We believe that a similar result holds for all random matrices having independent entries with mean zero and bounded variances.

\begin{conjecture} \label{conj:conj1} Let $c \le C$ be positive constants.
Let $A$ be the random matrix of size $n$ whose entries are independent random variables with
mean zero and variances between $c$ and $C$. Then with probability tending to one

\begin{equation}
|\det A| = n^{O(1)} \E |\det A |, \, \, \, \det A ^2=n^{O(1)}\E(\det A^2)
\end{equation}
\end{conjecture}

This conjecture looks highly  non-trivial. As a first step, we  consider the case
when the entries of $A$ are scaled Bernoulli random variables (namely, the ${ij}$ entry  takes values $\pm c_{ij}$ with probability half). Our experience is that this is usually the hardest case and its understanding would lead to the solution of the general case.  Our main result is

\begin{theorem} \label{thm:mainresult}
Let $0<c<C$ and $B>0$ be fixed.  Let $A$ be a random $n \times n$ matrix matrix whose entries $a_{ij}$ takes values $\pm c_{ij}$ with probability $1/2$, independently,  where  $ c\le |c_{ij}|  \le C$. Then with probability $1-n^{-B}$,
\begin{equation*}
|\det A|=\exp(O (n^{2/3} \log n)) \E(|\det A|),
\end{equation*}

\noindent and

\begin{equation*}
 \det (A^{2})= \exp(O(n^{2/3} \log n))\E(\det A^2).
\end{equation*}

Here the hidden constants in the $O$ notation may depend on $c, C, B$.

\end{theorem}

In the case $c=C=1$ (i.e., the entries of $A$ are iid Bernoulli), a better bound
$\exp( O(\sqrt {n \log n} )$ was recently proved in \cite{TVdet}. The approach in \cite{TVdet}, however, does not extend to random matrices with entries having different variances.
In the present approach,  it seems to require some new ideas  in order to significantly reduce the constant $2/3$.

If one assumes that the entries of $A$ are Gaussian (with different variances $c_{ij}^2$), then a weaker bound ($\exp(\epsilon n)$ for any positive $\epsilon$) was proved by Friedland, Rider and Zeitouni \cite{FRZ}. Our Theorem \ref{thm:mainresult} also holds for this case, with the same proof (see Section \ref{section:general}) and thus we obtain an improvement for the main result of \cite{FRZ}.

\section{Computing permanents}

Let us now consider
$\det M$ and $\per M$ from the computational point of view. It is not hard to compute $\det M$. In fact, there are effective algorithms to compute the whole spectra of $M$. The problem of computing $\per M$, on the other hand, is notoriously hard, and has been a challenge in theoretical computer science for many years.

A well-known observation that relates the problem of computing the permanent to that of determinant is the following. Let $u_{ij}$ be independent random variables with mean zero and variance one. Given a matrix $M$ with entries $a_{ij}$, define a random matrix $A$ with entries
$\sqrt {a_{ij}} u_{ij}$. Then, using linearity of expectation, it is easy to verify that
\begin{equation} \label{eqn:identity}
\E(\det A^2)=\textrm{per}(M).
\end{equation}

If $\det A^2$ is strongly concentrated around its mean, then
\eqref{eqn:identity} leads to the following very simple algorithm:
Given $M$, create a random sample of $A$. Compute $\det A^2$ and output it as an estimator for $\per A$. The core of the analysis is then to bound the degree of concentration of  $\det A^2$ around its expectation.

We mention here that in the case when $M$ has non-negative entries,
 the famous work of Jerrum and Sinclair \cite{JS} and Jerrum, Sinclair, Vigoda \cite{JSV}  gave an fully polynomial randomized approximation scheme for the problem, using the Markov-chain Monte Carlo approach. Theoretically, this result is as good as
it gets. On the other hand, the determinant estimator approach is still of interest, thanks to its
simplicity and implementability. (The  Markov chain algorithm requires running time $\Theta (n^7)$.)

In \cite{GG}, Godsil and Gutman proposed setting $u_{ij}$ to be iid Bernoulli random variables. Following the literature, we
call this algorithm the Godsil-Gutman estimator. This is perhaps the simplest
estimator. On the other hand, its analysis
 seems non-trivial. To illustrate this, let us consider the case when $M$ is the all-one matrix. Clearly $\per M= n!$. On the other hand, it is already not easy to prove that with high probability
$\det A \neq 0$ (this was first done by Koml\'os \cite{Kom1}). Effective bounds on $|\det A|$ have only recently become known(see \cite{TVdet}).


If one forces $u_{ij}$ to have a continuous distribution, the situation is more favorable. For instance, it is trivial that $\det A \neq 0$ with probability one. By setting $u_{ij}$ to be iid Gaussian variables, Barvinok \cite{Bar} showed that one can approximate the permanent of a non-negative matrix within a factor of $c^n$, for some constant $0<c<1$. A well-known problem with using Gaussian (or continuous) distribution is that in practice the implementation involves a truncated version of each variable. If the goal function (which is a function of many random variables) has a small Lipschitz coefficient, then this routine is effective. However, if its Lipschitz coefficient is large, then one needs to use a very fine approximation, and this increases the complexity of the input and would raise some challenges in implementation.

It is known that if one allows the matrix to have zero entries, then determinant estimators do not necessarily give a good approximation to the permanent.  For example,
Barvinok gave an example where the permanent is $2^n$ but the Godsil-Gutman estimator almost always returns $0$, and another where his own estimator will almost surely perform no better than an $\exp(O(n))$ approximation.
On the other hand,  Friedland, Rider and Zeitouni\cite{FRZ} showed that if the entries are strictly bounded from above and below by positive constants, then Barvinok  estimator gives an approximation factor  $\exp(\eps n))$, for any fixed $\eps >0$.

As a consequence of Theorem \ref{thm:mainresult} and Theorem \ref{thm:genresult}, we obtain
the following improvement

\begin{theorem} \label{theorem:appoximation}
Let $A$ be a (deterministic) square matrix of size $n$ with entries between $c$ and $C$, where $c$ and $C$ are positive constants. Then both the Godsil-Gutman and Barvinok estimators approximate $\per A$ within a  factor of $\exp( n^{2/3} \log n)$ with probability tending to one. \end{theorem}

If Conjecture \ref{conj:conj1} holds, then one can
improve the approximation factor to $n^{O(1)}$.

It remains  a tantalizing problem to analyze the determinant estimator for the case when the entries
of $A$ are not non-negative real numbers. Notice that \eqref{eqn:identity} still holds in this case, but no effective algorithm is known.

\section {The main ideas}

We start with the well-known identity

\begin{equation} \label{eqn:detprod}
 \det A^2=\det(AA^T)=\prod_{i=1}^n \sigma_i^2
\end{equation}
where $0 \leq \sigma_1 \leq \sigma_2 \leq \dots \leq \sigma_n$ are the singular values of $A$.

 If one could show that each singular value $\sigma_i$ is very strongly concentrated around some non-zero value, then
 $\det A^2$ would be so as well. Unfortunately, such a result is not available.
 In \cite{AKV}, it was shown, via Talagrand's inequality, that the largest singular values are  strongly concentrated, but
 the degree of concentration decreases rather quickly as the index decreases.

 To overcome this obstacle, we will follow the approach in \cite{FRZ}, which is based on the fact that, roughly speaking, the counting measure generated by the $\sigma_i$ is strongly concentrated.
 This fact was proved by Guionnet and Zeitouni in an earlier paper \cite{GZ}, also using Talagrand's inequality. Guionnet and Zeitouni's result asserts  that (after a proper normalization by a factor $1/\sqrt n$) any fixed interval, with high probability, contains the right number of singular values.  This enables one to show that
 the product of  most of the singular values is close to the expectation.

 The main technical barrier of this approach arises at the end of the spectrum.  The Guionnet-Zeitouni
 result does not reveal any information about the few smallest  singular values. In \cite{FRZ}, the authors needed to exploit  the Gaussian assumption (following an approach of Bai \cite{Bai}) in order to take care of these singular values. This technique, however, is not applicable for discrete distributions such as Bernoulli. In particular, it does not even show that
a  random matrix with discrete entries is non-singular with high probability.

 The proof of Theorem \ref{thm:mainresult} requires two new ingredients. The first is a lower bound on the smallest singular value $\sigma_n$. In \cite{TV2}, it was  shown, for many models of random matrices that  $\sigma_n$ is at least $n^{-C}$, for some constant $C$. While the models  in \cite{TV2} do not include the type of random matrices we consider here,
 we are able to  modify the proof, without too many difficulties, to treat our case.

  To continue, naturally one would try to use
   the uniform bound $n^{-C}$ for all  singular values which have not been treated by the concentration result. These will be singular values which are less by some threshold $\eps(n)$.
   It is now critical to estimate the number of such singular values.
   The value of $\eps(n)$ will be too small for the concentration result of Guionnet and Zeitouni to give information about this number.   The second main ingredient of our proof is a method that provides a good bound. This is based on a simple, but useful, identity (discovered in \cite{TVuniv}) which gives a relation between the singular values $\sigma_i$ and the distances $d_i$.

\section{A more general theorem and the main lemmas}
We will actually prove the following more general case of Theorem \ref{thm:mainresult}, where we merely require the entries to be bounded and have bounded variance instead of to be Bernoulli.

\begin{theorem} \label{thm:genresult}
 Let $K>0, B>0$, and $0<c<1$ be fixed.  Let $A$ be a random
  $n \times n$ matrix whose entries $a_{ij}$ are random variables
  satisfying
\begin{itemize}
 \item{$c \leq \V(a_{ij}) \leq \frac{1}{c}$}
 \item{$\P(|a_{ij} - E(a_{ij})| \leq K)=1$.}
\end{itemize}
Then  with probability $1-n^{-B}$,
\begin{equation*}
|\det A| =\exp(O(n^{2/3} \log n))\E(|\det A|)
\end{equation*}
and
\begin{equation*}
\det A^2 =\exp(O(n^{2/3} \log n))\E(\det A^2),
\end{equation*}

where the constant implicit in the $O$ notation depends on $K, B,$ and $c$.
\end{theorem}
\begin{remark}
 Here and elsewhere the relation $a_n=O(b_n)$ indicates that the ratio $a_n/b_n$ is bounded above in absolute value as $n$ tends to infinity.  In particular, the theorem above gives both an upper bound and a lower bound on the ratio between the determinant and its expectation.
\end{remark}

\begin{remark}
 The uniform boundedness condition can be replaced by the condition that all of the entries have a Gaussian distribution; see section 8.
\end{remark}

Recall \eqref{eqn:detprod},
\begin{equation}
 (\det A)^2=\det(AA^T)=\prod_{\sigma \in \spec{AA^{T}}}= \prod_{i=1}^n \sigma_i^2
\end{equation}
where $0 \leq \sigma_1 \leq \sigma_2 \leq \dots \leq \sigma_n$ are the singular values of $A$.

We will start in a similar way as in \cite{FRZ}. Let
$\epsilon$ be a parameter to be determined later (which may depend on $n$). We estimate
\eqref{eqn:detprod} by dividing the spectrum into two parts, writing
$|\det A|=\det_{trunc} A \det_{small} A$, where

\begin{eqnarray*}
 \det_{trunc}&=&\left(\prod_{\sigma \in \spec(AA^T)} \max \{\sigma, \epsilon^2\}\right)^{1/2},  \\
\det_{small}&=&\left(\prod_{\sigma\in \spec(AA^T)} \min\{\sigma \epsilon^{-2}, 1\}\right)^{1/2}.
\end{eqnarray*}

We show that  $ \det_{trunc} A$ and $ \det_{trunc} A^2$ are strongly concentrated around their means

\begin{lemma} \label{lemma:dettrunc}
 There is a constant $c_0>0$ dependent only on $c$ such that
\begin{equation*}
 \det_{trunc} A= \exp(O(n  \epsilon^{-2} \log n)) \E(\det_{trunc} A)
\end{equation*}
and
\begin{equation*}
 \det_{trunc} A^2= \exp(O(n  \epsilon^{-2} \log n)) \E(\det_{trunc} A^2)
\end{equation*}
with probability $1-O(n^{-c_0 \log n})$.
\end{lemma}

The proof is presented in Section 3.

To handle $\det_{small} A$, notice that
\begin{equation} \label{eqn:detsmall}
 1 \geq \det_{small} A \geq \min\{1, (\sigma_n(A) \epsilon^{-1})^{s_{\epsilon}(A)}\}
\end{equation}
where $\sigma_n(A)$ is the smallest singular value of $A$ and $s_{\epsilon}(A)$ denotes the number of singular values of $A$ which are at most $\epsilon$.  We can therefore bound $\det_{small} A$ from below by using the following two lemmas.

\begin{lemma} \label{lemma:smallestsv}
 For any $B>0$,
\begin{equation*}
 \P(\sigma_n(A)<n^{-4B-7}) \leq n^{-B}
\end{equation*}

\end{lemma}
We remark that  $-4B-7$ is pretty far from being optimal and can be improved, but doing so would not affect our final results in any essential way.
This lemma  is a variant of many results proved in \cite{TVpert} (see also \cite{TVcirc}).   However, \cite{TVpert} required that the distributions of the entries of $A$ to be dominated in a certain Fourier analytic sense by a single common distribution.  Our matrices do not satisfy this assumption. However, we are able to modify the proof, without too many difficulties, to obtain the desired result.

\begin{lemma}\label{lemma:fewsmallsv}
 Let $r \geq \log^4 n$, and assume $c \leq \V(a_{ij}) \leq \frac{1}{c}$.  Then
\begin{equation*}
\P\left( \sigma_{2r} (A) \leq \frac{rc^2}{2\sqrt{n-r}} \right)=o(n^{-\log n}),
\end{equation*}
\end{lemma}

The above two lemmas combine to show that no singular value of $A$ is likely to be so small as to have too large an effect on the determinant, and, furthermore, we can also deduce that not many singular values will have to be handled by $\det_{small}$.

Let us for now assume the previous two lemmas to be true.  By taking $r=\frac{3\epsilon \sqrt{n}}{c^2}$ in Lemma \ref{lemma:fewsmallsv}, we see that with high probability $s_{\epsilon}(A)=O(n^{1/2} \epsilon)$.  Combining this with Lemma \ref{lemma:fewsmallsv} and the bounds in \eqref{eqn:detsmall}, we see that for any $B>0$ we have with probability $1-n^{-B+o(1)}$ that
\begin{equation*}
\det_{small} A=\exp(O(n^{1/2} \epsilon \log n )),
\end{equation*}
which therefore implies that with the same probability
\begin{equation} \label{eqn:detineq}
 \det_{trunc} A \geq |\det A| \geq \exp(O(n^{1/2} \epsilon \log n)) \det_{trunc} A.
\end{equation}
(note again that the use of $O$ in the lower bound here indicates an exponent bounded in magnitude.) Now let us
fix $\epsilon=n^{1/6}$. Combining the second half of the above
inequality with Lemma \ref{lemma:dettrunc}, we see that with
probability $1-n^{-B+o(1)}$ we have
\begin{equation*}
 |\det A| \geq \exp(O(n^{2/3} \log n)) \E(\det_{trunc} A).
\end{equation*}
Taking expectations, we find that
\begin{equation} \label{eqn:expdetineq}
 \E(\det_{trunc} A) \geq \E(|\det A|) \geq (1+o(1)) \exp(O(n^{2/3} \log n)) \E(\det_{trunc} A).
\end{equation}
The first half of Theorem \ref{thm:genresult} follows from combining \eqref{eqn:detineq}, \eqref{eqn:expdetineq}, and Lemma \ref{lemma:dettrunc}.  The second half follows from the identical argument being applied to $\det A^2$.

\section{The Proof of Lemma \ref{lemma:dettrunc}}
As in \cite{FRZ}, we begin with the spectral concentration results of Guionnet and Zeitouni, in particular the following special case of Corollary 1.8(a) in \cite{GZ}:
\begin{theorem}
 Let $Y$ be an $n \times n$ matrix whose entries are independent random variables each having support on a compact set of diameter at most $K$, and let $Z=Y^T Y$.  Let $\lambda_1, \dots, \lambda_n$ be the eigenvalues of $Z$, and let $f$ be an increasing, convex, function such that $g(x)=f(x^2)$ has Lipschitz norm $|g|_L$.  Then for any $\delta>\delta_0:=\frac{2K\sqrt{\pi}|g|_L}{n}$,
\begin{equation*}
 \P(|\sum_{i=1}^n f(\lambda_i) -\E(\sum_{i=1}^n f(\lambda_i))| > 2\delta n) \leq 4\exp\left(-\frac{(\delta-\delta_0)^2 n}{K^2 |g|_L^2}\right).
\end{equation*}
 \end{theorem}

Ideally, we would like to apply this theorem with $f$ taken to be the logarithm, so that $\sum f(\lambda_i)=\log \det A^2$.  The difficulty is that the logarithm is not Lipschitz.  To overcome this problem, we
follow \cite{FRZ} and truncate the logarithm.  Write
\begin{equation*}
 \log^{\epsilon}x=\max\{2 \log \epsilon, \log x\},
\end{equation*}
where $\log^{\epsilon}(0)$ is defined to be $2 \log \epsilon$.  Note that we have
\begin{equation} \label{eqn:spec}
 \log(\det_{trunc} A )=\frac{1}{2} \sum_{\sigma \in \spec(AA^T)} \log^{\epsilon}(\sigma).
\end{equation}

Although $\log^{\epsilon}(x^2)$ now has finite Lipschitz constant $\frac{1}{\epsilon}$ (this was the purpose of truncating the logarithm), it is not convex.  However, it can easily be written as the difference of two convex Lipschitz functions, so the above theorem applies, and we have for some absolute constants $C_0$ and $C_1$ and any $\delta>=\delta_0:=C_0 \epsilon^{-1}/n$ that


\begin{equation} \label{eqn:logbound}
 \P(|\log(\det_{trunc} A) - \E(\log(\det_{trunc} A))| > \delta n) \leq 4 \exp(-C_1\frac{\epsilon^2 n \delta^2c^2}{16})
\end{equation}
 Taking $\delta=\frac{\log n}{\epsilon \sqrt{n}}$, we see that for some constant $c_0$ we have 
\begin{equation} \label{eqn:dettruncext} 
 \P(|\log(\det_{trunc} A)-\E(\log(\det_{trunc} A))| > \frac{\sqrt{n} \log n }{\epsilon} =O(n^{-c_0 \log n})
\end{equation}

This would be exactly the result we wanted, if only the expectation and the logarithm were switched on the left hand side of \eqref{eqn:logbound}.  Following \cite{FRZ}, we now write
\begin{equation*}
 U(A)=\log(\det_{trunc} A)-\E(\log(\det_{trunc} A)).
\end{equation*}
 We know $\E(U)=0$, and by Jensen's inequality we have
\begin{equation*}
1 \leq \E(e^U) \leq \E(e^{|U|}) \leq 1+\int_0^{\infty} e^t \P(|U|>t) dt.
\end{equation*}

It follows from the above and \eqref{eqn:logbound} that
\begin{equation*}
 1 \leq \E(e^U) = \frac{\E(\det_{trunc} A)}{e^{\E(\log(\det_{trunc} A))}} \leq exp(O(\frac{n}{\epsilon^2}))
\end{equation*}
and the first half of Lemma \ref{lemma:dettrunc} follows by taking logarithms and combining with \eqref{eqn:dettruncext}.  The second half follows from the identical calculation applied to $e^{2U}$.

\begin{remark}
If we only had required that the truncated determinant concentrate \textit{somewhere}, the argument above would have given a stronger bound (roughly $\exp(\epsilon^{-1} n^{1/2} \log n )$).  The dominant term in our bound came from showing that the ``somewhere'' was close to the actual expectation.

Also, we did not at any point use our lower bound on the variance of the entries.  In particular, this truncated determinant will be concentrated around its expectation even if we allow most of the entries of $A$ to be non-random.  However, it is not true in general that $\det_{trunc} A$ will be close to $\det A$.
\end{remark}
\section{The Proof of Lemma \ref{lemma:smallestsv}}
We begin by first reducing from the general case back to the case of
Bernoulli Matrices. To do so, we will use the idea of Bernoulli
decomposition from a paper of Aizenman et. al. \cite{AFKW}.  In this
paper, it is shown that for any random variable $X$ that is
nondegenerate (not taking on any single value with probability 1),
we can find a $p \in (0,1)$ and functions $f(t)$ and $g(t)$ such
that
\begin{itemize}
\item{If $t$ is uniform on $[0,1]$, and $\epsilon$ is a Bernoulli variable independently equal to 1 (with probability $p$) or 0 (with probability $1-p$), then $f(t)+g(t)\epsilon$ has the same distribution as $X$}
\item{$\inf g(t)>0$.}
\end{itemize}
Recall that we are assuming that our entries are both uniformly bounded in magnitude by $K$ and bounded below in variance by $c$.  It follows from the methods of \cite{AFKW} (see Remark 2.1(i) there), that in this case we can find a Bernoulli composition $a_{ij}=f_{ij}(t_{ij})+g(t_{ij})\epsilon_{ij}$ of every entry of $A$ in which the $g(t_{ij})$ have a uniform lower bound $\beta=\beta(K,c)$ for all values of $i, j,$ and $t$, and for which the $p_{ij}$ in the decompositions are uniformly bounded away from $0$ and $1$.

We now view our matrix as being formed in two steps.  First, we expose $t_{ij}$ for each entry.  At this point every entry can be viewed as having a shifted Bernoulli distribution.   Next we expose the $\epsilon_{ij}$.  It follows by taking expectations over all possible values of $t_{ij}$ that it suffices to show the following
\begin{lemma} \label{lemma:Bernresult}
 Let $0<q<\frac{1}{2}$ and $B,C,c>0$ be fixed.  Let $A$ be a matrix whose entries are independent random variables distributed as $a_{ij}=m_{ij}+\epsilon_{ij} n_{ij}$, where $|m_{ij}|<n^{1/8}$ and $c<n_{ij}<C$, and furthermore the $\epsilon_{ij}$ satisfy
\begin{equation*}
 q<\P(\epsilon_{ij}=1)=1-\P(\epsilon_{ij}=-1)<1-q.
\end{equation*}
Then for sufficiently large $n$ we have
\begin{equation*}
 \P(\sigma_n(A)<n^{-4B-7}) \leq n^{-B}.
\end{equation*}
\end{lemma}
\begin{remark}
The form of this theorem is very similar to that of the smoothed analysis of the smallest singular value in \cite{TVpert}.  The key difference here is that we no longer require the $n_{ij}$ to be identical.
\end{remark}

Proving Lemma \ref{lemma:Bernresult} is equivalent to bounding the probability that for some unit vector $v$ we have $||Av|| \leq n^{-4B-7}$.  We will do this by dividing the vectors into two classes, which should be thought of as ``structured'' and ``unstructured'', for an appropriate definition of ``structured'' depending both on $A$ and on $B$.

\begin{definition}
A vector $v$ is \textbf{rich} if there is some $i$  for which
\begin{equation*}
\sup_z \P(|\sum_{j=1}^n a_{ij} v_j-z|<n^{-4B-13/2}) \geq n^{-B-1}
\end{equation*}
Otherwise $v$ is \textbf{poor}.
\end{definition}

Equivalently, a poor vector is one for which no individual coordinate of $Av$ is too concentrated.  Lemma \ref{lemma:smallestsv} would be an immediate consequence of the following two lemmas.
\begin{lemma} \label{lemma:poorv}
\begin{equation*}
\P(||Av|| \leq n^{-4B-7} \textrm{ for some poor } v)\leq \frac{1}{2}n^{-B}
\end{equation*}
\end{lemma}
\begin{lemma} \label{lemma:richv}
\begin{equation*}
\P(||Av|| \leq n^{-4B-7} \textrm{ for some rich } v)\leq \frac{1}{2}n^{-B}
\end{equation*}
\end{lemma}

\textbf{ Proof of Lemma \ref{lemma:poorv}}:

We adapt an argument from  \cite{Rud} (see also \cite{ TVpert}).  Let $E$ be the event that for some poor unit vector $v$ we have $||Av|| \leq n^{-4B-7}$.  If $E$ holds, then the least singular value of $A$ is at most $n^{-4B-7}$, so the same must hold for $A^T$.  For $1 \leq j \leq n$, let $F_j$ be the event that there exists a unit vector $w=(w_1, \dots w_n)^T$ which simultaneously satisfies
\begin{equation*}
 ||w^T A|| \leq n^{-4B-7}, \, \, \, |w_j| \geq \frac{1}{\sqrt{n}}.
\end{equation*}
Since every $w$ has at least one coordinate at least $n^{-1/2}$ in magnitude, we have
\begin{equation*}
 \P(E) \leq \sum_{i=1}^n \P(E \wedge F_j).
\end{equation*}
Now let $j$ be fixed.  Let $A_1, \dots A_n$ be the rows of $A$. We will condition on all of the rows except row $j$.  If $E$ is to hold, there must be a poor $v$ such that
\begin{equation*}
 (\sum_{i=1}^n |A_i \cdot v|^2)^{1/2} = || A v|| \leq n^{-4B-7}.
\end{equation*}
It follows that if $\P(E|A_1, \dots A_{j-1}, A_{j+1} ,\dots A_n)$ is non-zero, then there is a poor $u$ such that
\begin{equation} \label{eqn:pooru}
 (\sum_{i \neq j} |A_i \cdot u|^2)^{1/2} \leq n^{-4B-7}
\end{equation}
Conversely, by our assumptions on $w$ we have that if $F_j$ holds, then
\begin{equation*}
 ||\sum_{i \neq j} w_i A_i || \leq n^{-4B-7}.
\end{equation*}
Taking inner products with $u$ and using the triangle inequality, we conclude
\begin{equation*}
 |w_j| |A_j \cdot u| \leq \sum_{i \neq j} |w_i| |A_i \cdot u|+n^{-4B-7}
\end{equation*}
Combining the above with \eqref{eqn:pooru}, the Cauchy-Schwartz inequality, and our assumption on $|w_j|$, we obtain that if both $E$ and $F_j$ hold, then
\begin{equation*}
 |A_j \cdot u| \leq 2n^{-4B+13/2}.
\end{equation*}
On the other hand, since $u$ is poor and $A_j$ and $u$ are independent, we have that
\begin{equation*}
 \P(|A_j \cdot u| \leq 2n^{-4B+13/2} | A_1, \dots A_{j-1}, A_{j+1}, \dots A_n) \leq n^{-B-1}.
\end{equation*}
Combining the above, we see that
\begin{equation*}
 \P(E \wedge F_j | A_1 \dots A_{j-1}, A_{j+1}, \dots A_n) \leq n^{-B-1},
\end{equation*}
regardless of our choice of the remaining $n-1$ rows.  It follows that $\P(E \wedge F_j) \leq n^{-B-1}$ for every $j$, and therefore that $\P(E) \leq n^{-B}$.

\textbf{ Proof of Lemma \ref{lemma:richv}}:
Let $J$ be the unique integer satisfying $2B+2<J \leq 2B+3$ and let $\delta=(B+1)/J$.  Let $\gamma>0$ be a constant chosen to be sufficiently small that $\delta+3\gamma<1/2$ and $(4B+7)\gamma<\frac{1}{2}$.  Finally, we let $D=2+2\gamma$.

Let $v$ be a rich unit vector.  We define
\begin{equation*}
 g(j):=\sup_{i, z} \P(|A_i^T v-z|< n^{-4B-13/2+Dj})
\end{equation*}
Clearly $0 \leq g(j) \leq 1$, and $g(j)$ is an increasing function in $j$.  The assumption that $v$ is rich is equivalent to the statement that $g(0) \geq n^{-B-1}$.  It follows from the pigeonhole principle that for some $0<j \leq J-1$ we have
\begin{equation*}
 g(j+1) \leq n^{\delta} g(j).
\end{equation*}
For $0 \leq j \leq J-1$ and $1 \leq k \leq \lceil \frac{(A+1)}{\gamma} \rceil$, we define $\Omega_{j,k}$ to be the collection of rich $v$ satisfying both 
\begin{equation*}
 g(j+1) \leq n^{\delta} g(j) \, \, \, \textrm{ and } \, \, \, g(j) \in [n^{-k\gamma},n^{-(k-1)\gamma}]
\end{equation*}

Since every rich $v$ is contained in some $\Omega_{j,k}$, and there are only a bounded number of pairs $(j,k)$, it suffices to prove that for every fixed $j$ and $k$ we have
\begin{equation} \label{eqn:omegajk}
 \P(||A v|| \leq n^{-4B-7} \textrm{for some }v \in \Omega_{j,k})=o(n^{-B}).
\end{equation}

Our goal will now be to construct a $\beta-$net for each $\Omega_{j,k}$, that is a set $V_0$ such that any point in $\Omega_{j,k}$ is within (Euclidean) distance $\beta$ of some point in $V_0$.  Assuming that for sufficiently small $\beta$ the net is not too large, we will then be able to obtain $\ref{lemma:richv}$ by a union bound.  We begin bounding the size of the net with the following result, a special case of \cite[Thm. 3.2, see also Remark 2.8]{TVcirc}:

\begin{theorem} \label{thm:lonetcirc}
 Let $0<q<\frac{1}{2}$ and let $x_1, \dots x_n$ be independent random variables taking on values in $\{1, -1\}$ and satisfying
\begin{equation*}
 q \leq \P(x_{j}=1) \leq 1-q.
\end{equation*}
Let $0<\delta<1$ be fixed, and let $p$ and $\beta$ be chosen to satisfy $p=n^{-O(1)}$ and $\beta>\exp(-n^{-\delta/2})$.  Then the set of vectors $(v_1, \dots v_n)$ satisfying
\begin{equation*}
 \sup_{z \in \C} \P(|\sum_{i=1}^n v_i x_i - z| < \beta)<p
\end{equation*}
has a $\beta$-net in the $l_{\infty}$ norm of size at most $n^{-(1/2+\delta)n}p^{-n} +\exp(o(n))$.
\end{theorem}

For any particular $i$, we have
\begin{equation} \label{eqn:betanet}
 \P(|\sum_{j=1}^n (m_{ij}+n_{ij} \epsilon_{ij}) v_j -z|<\beta) = \P(|\sum_{j=1}^n \epsilon_{ij} \widetilde{v_j} -\widetilde{z}|<\beta,
\end{equation}
where $\widetilde{v_j}=v_j n_{ij}$ and $\widetilde{z}=z-\sum_j m_{ij} v_j$.  For any particular coordinate of $Av$, Theorem \ref{thm:lonetcirc} gives an upper bound on the minimal size of a $\beta-$net for the set of $\widetilde{v}$ for which the right hand side of \eqref{eqn:betanet} holds with probability at least $p$.  By taking an affine transformation $v_j \rightarrow \frac{v_j}{n_{ij}}$ of the case $\beta=n^{-4B-13/2+Dj}, p=n^{-k\gamma}$ of this net and taking the union of the resulting net for each coordinate, we obtain the following modified version of Theorem \ref{thm:lonetcirc}:

\begin{lemma} \label{lemma:lonetadjusted}
 Let $x_1, \dots x_n$ be independent and have the form $x_i=m_i+\epsilon_i n_i$, where the $m,n,\epsilon$ are as in Lemma \ref{lemma:Bernresult}.  Let $0<\delta<1$ be fixed. Then $\Omega_{j,k}$ has an $\frac{n^{-4B-13/2+Dj}}{c}$-net in the $l_{\infty}$ norm of size at most $n^{1-(1/2+\delta)n} n^{k\gamma n} +\exp(o(n))$.
\end{lemma}

Let $V_0$ be a net guaranteed by the above lemma, and consider any $v' \in V_0$ and $v \in \Omega_{j,k}$ such that $||v-v'||_{\infty} \leq \beta$.  Our bounds on the $m_{ij}$ and $n_{ij}$ guarantee that (assuming $n$ to be sufficiently large) the spectral norm of $A$ satisfies $\sigma_n(A)<n$.  Since
\begin{equation*}
 ||Av'|| \leq ||Av||+||A(v-v')|| \leq ||Av||+n^{1/2} \sigma_n(A) ||v-v'||_{\infty},
\end{equation*}
it follows that if $||Av|| \leq n^{-4B-7}$ then
\begin{equation*}
 ||Av'|| \leq (1+\frac{1}{c}) n^{-4B-4+Dj}.
\end{equation*}

It follows that there must be at least $n-n^{1-\gamma}$ rows of $A$ for which
\begin{equation} \label{eqn:smallrow}
 |A_i^T v'| \leq n^{-4B-9/2+Dj+\gamma}.
\end{equation}

On the other hand, we also have for any $i$ for which \eqref{eqn:smallrow} holds that
\begin{eqnarray*}
 |A_i^Tv| &\leq& |X_i^T v'|+||v-v'||_{\infty} \sum_{j=1}^n (m_{ij}+n_{ij}) \\
&\leq& n^{-4B-9/2+Dj+\gamma}+\frac{1}{c}n^{-4B-9/2+Dj}(\frac{1}{c}+n^{1/8}) \\
&\leq& n^{-4B-7+D(j+1)}.
\end{eqnarray*}
Where the last inequality comes from our definition of $D$.

It follows that
\begin{eqnarray*}
 \P(|A_i^Tv'| \leq n^{-4B-9/2+Dj+\gamma}) &\leq& \P(|A_i^Tv| \leq n^{-4B-7+D(j+1)}) \\
&\leq& n^{\delta} g(j) \\
&\leq& n^{\delta+\gamma-k\epsilon}
\end{eqnarray*}
where for the last two inequalities we use the definition of $\Omega_{j,k}$.

This will be sufficient to handle the case where $k$ is sufficiently large.  For smaller $k$, we note that by our choice of $\gamma$ and $D$ we have $-4B-9/2+Dj+\gamma<-1$, so
\begin{equation*}
 \P(|A_i^Tv'| \leq n^{-4B-9/2+Dj+\gamma}) < \P(|A_i^T v' < \frac{1}{n}|),
\end{equation*}
which can easily be checked to be at most $1-q$.

Therefore
\begin{equation*}
 \P(|A_i^T v'| \leq n^{-4B-9/2+Dj+\gamma)}) \leq \min(n^{\delta+\gamma-k\epsilon}, 1-q).
\end{equation*}
Taking the union bound over all sets of $n-n^{1-\gamma}$ rows, we see that
\begin{equation*}
 \P(||Av'|| \leq n^{-4B-4+Dj+\gamma}) \leq \min(n^{\delta+\gamma-k\epsilon}, 1-q)^{n'} \binom{n}{n-n^{1-\gamma}}
\end{equation*}
for any particular $v'$ in our net for $\Omega_{j,k}$.  Taking the union bound over the entire net, we obtain that the probability that the left hand side of \eqref{eqn:omegajk} holds is at most
\begin{equation*}
 (n^{-(1/2+\delta)n+1} n^{k\gamma n} +\exp(o(n))) \min(n^{\delta+\gamma-k\epsilon}, 1-q)^{n-n^{1-\gamma}} \binom{n}{n-n^{1-\gamma}}
\end{equation*}
which can be verified to be exponentially small by a routine calculation.

\section{The Proof of Lemma \ref{lemma:fewsmallsv}}
As in the proof of Lemma \ref{lemma:smallestsv}, it suffices by Bernoulli decomposition to prove the following special case of this lemma:

\begin{lemma} \label{lemma:Bernresult2}
 Let $0<q<\frac{1}{2}$ and $B,C,c>0$ be fixed. Let $A$ be a matrix whose entries are independent random variables having distributed as $a_{ij}=m_{ij}+\epsilon_{ij} n_{ij}$, where $|m_{ij}|<n^{1/8}$ and $c<n_{ij}<C$, and furthermore the $\epsilon_{ij}$ satisfy
\begin{equation*}
 q<\P(\epsilon_{ij}=1)=1-\P(\epsilon_{ij}=-1)<1-q.
\end{equation*}
Then for $r \geq \log^4 n$,
\begin{equation*}
 \P\left(\sigma_{2r}(A)<\frac{rc^2}{2\sqrt{n-r}}\right) =o(n^{-\log n}).
\end{equation*}
\end{lemma}

To prove this Lemma we are going to use the following lemma from \cite{TVuniv}

\begin{lemma} \label{lemma:identity} \cite{TVuniv} Let $M$ be an $m \times n$ matrix ($m \le n$). Let $d_i$ be the distance from its $i$th row vector to the space spanned by the first $i-1$ rows and $\sigma_i$ be its singular values. Then
$$\sum_{i=1}^m d_i^{-2} =\sum_{i=1}^m \sigma_i^{-2} $$

\end{lemma}
Recall that $\log^4 n<r<\frac{n}{2}$.  By the interlacing inequalities for singular values (see, for example,  Theorem 7.3.9 in \cite{hj}), we have that
\begin{equation*}
\sigma_{2r}(A) \geq \sigma_r (A'),
\end{equation*}
where $A'$ is the matrix formed by removing the last $r$ columns from $A$.

To bound the right hand side of this equation, we note that
\begin{eqnarray}
 \sigma_r (A')^{-2} &\leq&  \frac{1}{r} \sum_{k=1}^r \sigma_k (A')^{-2} \nonumber \\
&\leq& \frac{1}{r}\sum_{k=1}^{n-r} \sigma_k(A')^{-2} \nonumber \\
&=& \frac{1}{r} \sum_{i=1}^{n-r} d_i^{-2}, \label{eqn:singjbound}
\end{eqnarray}
where $d_i$ denotes the distance from the $i^{th}$ column of $A'$ to the span of the remaining columns, and for the last equality we use Lemma \ref{lemma:identity}.

Informally, this states that if a matrix has many small singular values, it must have many columns which are very close to the subspace spanned by the other columns.  Since $r$ is becoming increasingly large, the co-dimension of this subspace is increasing as well, so this should become unlikely.

Now let $i$ be fixed.  To bound the probability that $d_i$ is small, we first expose the subspace $S_i$ spanned by the remaining $n-r-1$ columns, then finally the remaining column.  Let $P$ denote the projection matrix onto that subspace, and let $p_{ij}$ be the entries of $P$.  Let $X_i=(a_{i1}, \dots a_{in})$ be this final column.  We have
\begin{eqnarray*}
 \E(d_i^2 | S_i) &=& \E(|X_i|^2) - \E(|P X_i|^2) \\
&=&\sum_{j=1}^n a_{ij}^2 - \sum_{k=1}^n \sum_{l=1}^n a_{ik} a_{il} p_{kl} \\
&=&\sum_{j=1}^n a_{ij}^2 (1-p_{kl}) \\
&\geq& \sum_{j=1}^n c^2 (1-p_{kl}) \\
&=& c^2 (n-Tr(P))=c^2 r
\end{eqnarray*}
as the terms with $k \neq l$ cancel by the independence of the entries of $A$, and we again use how the entries of $A$ are bounded away from 0.
\begin{remark}
 If the entries of $A$ were to have equal variance $c$, then the inequality here would actually be an equality, and the expected square distance would be independent of $S$.  This is what enables the arguments of \cite{Gir, TVdet} (which are based on row by row exposure of the matrix in question), and why those arguments don't carry over here to give an immediate estimate on the determinant of $A$.
\end{remark}

In other words, a random vector from our distribution will on average be far away from any fixed $n-r-1$ dimensional subspace.  It remains to show that it will typically be far away.

It follows from Talagrand's inequality \cite{Tal} and our bounds on the $n_{ij}$ that if $M_{i,S}$ is the median value of $d_i$ conditioned on $S_i$, then
\begin{equation*}
 \P(|d_i - M_{i,S}| \geq t | S_i) \leq 4\exp(-\frac{t^2}{64C^2}).
\end{equation*}
By an argument identical to that in \cite{TVdet}, it can be shown that $|M_i-\sqrt{E(d_i^2)}|\leq \frac{C^2}{c^2}$.  It therefore follows that for sufficiently large $r$,
\begin{equation*}
 \P(d_i \leq \frac{c^2 \sqrt{r}}{2} | S_i) \leq 4 \exp(-\frac{r^2 c^2}{300 C^2})=o(n^{-\log n-1}),
\end{equation*}
as by assumption $r>\log(n)^2$.

In particular, with probability $1-o(n^{-\log n})$ every $d_i$ will be at least $c^2 r/2$.  Combining this with \eqref{eqn:singjbound}, we see
\begin{eqnarray*}
 \P(\sigma_{2r} (A) \leq \frac{r c^2}{2 \sqrt{n-r}}) &\leq& \P(\sigma_r(A')^{-2} \geq \frac{4(n-r)}{c^4 r^2}) \\
&\leq& \P(\sum_{i=1}^{n-r} d_i^{-2} \geq \frac{4(n-r)}{c^4 r^2} \\
&=&o(n^{-\log n})
\end{eqnarray*}

\section{Concentration of Determinants for Gaussian variables} \label{section:general}

Although we have focused mainly on the concentration of determinants for the case of matrices with uniformly bounded entries, our main results also hold in the case where every entry has a Gaussian distribution, assuming that the means of the entries are uniformly bounded and the variances of the entries are uniformly bounded above and below.  In particular, this implies that our bounds hold for Barvinok's as well as Godsil and Gutman's estimator for the permanent.

The proof of Lemma \ref{lemma:dettrunc} is exactly the same as before, except that we use Corollary 1.8b of \cite{GZ} instead of Corollary 1.8a.  For the remaining two lemmas, we again use the idea of Bernoulli decomposition.  It can be explicitly checked that if $X$ is a Gaussian variable satisfying $|\E (X)|<K$, $c<\V(X)<C$, then $X$ can be decomposed as
\begin{equation}
 X=f(t)+g(t) \epsilon,
\end{equation}
where $t$ is uniform on $[0,1]$, $\epsilon$ is uniform on $\{-1, 1\}$.  Furthermore, we can do so in such a way that $g(t)$ is bounded uniformly from below, and the measure of the set of $t$ for which $|f(t)|+|g(t)|<\log^2 n$ is $o(n^{-\log n})$.

We now expose $t_{ij}$ for each entry of $A$.  At this point every entry will be a Bernoulli distribution and (except for an exceptional set of probability $o(n^{-B})$ for any $B$) the mean and variance of the entries will be bounded by $\log^2 n$.  Lemmas \ref{lemma:smallestsv} and \ref{lemma:fewsmallsv} now follow as before from Lemmas \ref{lemma:Bernresult} and \ref{lemma:Bernresult2}.

\textbf{Acknowledgements:} The authors with to thank the anonymous referees for their careful reading of and helpful comments on this paper.

\appendix
\section{Log-Normality of the Determinant of Gaussian Matrices}
In this appendix we will show that the determinant of a matrix whose entries are iid standard Gaussian random variables has the distribution given by Theorem \ref{theorem:normal}.  Our starting point is the formula
\begin{equation} \label{eqn:distprod}
 |\det A|=\prod_{i=1}^n d_i,
\end{equation}
where $d_i$ is the distance from the $i^{th}$ row of $A$ to the subspace spanned by the previous $i-1$ rows.

This formula is particularly useful for Gaussian vectors due to their rotational invariance: If $x$ is a random vector whose coordinates are iid Gaussian variables having mean zero, then distribution of the distance from $x$ to a fixed subspace $S$ is dependent only on the dimension of $S$.  If the dimension is $n-k$, then the distribution of the square of the distance follows a chi-square distribution with $k$ degrees of freedom.  In particular, this implies that the distribution of the determinant is the same as that where we treat each of the variables in \eqref{eqn:distprod} as being independent and following a chi distribution.  We will do so for the remainder of this appendix.

Taking logarithms in \eqref{eqn:distprod} and rearranging, we see that
\begin{equation*}
 2 \log (|\det A|)-\log n!=\sum_{i=1}^n \log(\frac{d_i^2}{i})=\sum_{i=1}^n \log(1+\frac{d_i^2-i}{i}).
\end{equation*}

Following the ideas of \cite{Gir}, we next perform a Taylor expansion on the right hand side, writing
\begin{equation} \label{eqn:CLTtaylor}
 \frac{2 \log (|\det A|)- \log n!}{\sqrt{2 \log n}}=\frac{\sum_{i=1}^n \frac{d_i^2-i}{i}}{\sqrt{2 \log n}} -\frac{1}{2} \frac{\sum_{i=1}^n (\frac{d_i^2-i}{i})^2}{\sqrt{2 \log n}}+\frac{1}{3} \frac{\sum_{i=1}^n (\frac{d_i^2-i}{i})^3}{\sqrt{2 \log n}}+\frac{\sum_{i=1}^n \epsilon_i}{\sqrt{2 \log n}}.
\end{equation}
We examine the terms in order.

It follows from standard facts about the chi-square distribution that $\frac{d_i^2-i}{i}$ has mean 0, variance $\frac{2}{i}$, and fourth moment $\frac{12i+48}{i^3}$.  It follows immediately from Lyapunov's Central Limit Theorem that the first term on the right hand side of \eqref{eqn:CLTtaylor} converges weakly to $N(0,1)$.

For the second term, we observe $(\frac{d_i^2-i}{i})^2$ has mean $\frac{2}{i}$ and variance $\frac{12i+48}{i^3}$.  In particular, the variance of $\sum (\frac{d_i^2-i}{i})^2$ is $o(\log n)$.  It follows that the second term converges to $\frac{-\log n}{\sqrt{2 \log n}}$.  Similarly, it follows from the moments of the chi-square distribution that the expectation and the variance of $\sum(\frac{d_i^2-i}{i})^3$ are $O(1)=o(\log(n))$, so the third term converges weakly to zero.

The final term is slightly more complicated due to the singularity of the logarithm at 0.  We first split the error term as $\epsilon_i=\epsilon_i^{'}+\epsilon_i^{''}$, where $\epsilon_i^{'}$ is zero whenever $d_i<\frac{i}{3}$, and $\epsilon_i^{''}$ is zero whenever $d_i \geq \frac{i}{3}$.  We will show separately that each of the contribution of each of these errors (converges weakly to zero).

For the first error term (the case where $d_i$ is large), we note that $|\epsilon_i^{'}|=O(|d_i^2-i|)^4$ and therefore (using the fourth moment given above) $\E(|\epsilon_i^{'}|)=O(1)$.  It follows that $\frac{\sum \epsilon_i^{'}}{\sqrt{2 \log n}}$ converges to zero, so the first part of our decomposition is negligible.

It can be checked by direct computation that $\E(|\log(d_i^2)|)$ is finite for any $i$.  The same therefore also holds for $\E(|\epsilon_i^{''}|)$, so it follows that for some function $s=s(n)$ diverging to infinity sufficiently slowly we have
\begin{equation} \label{eqn:gausstrunc}
 \frac{\sum_{i=1}^s \epsilon_i^{''}}{\sqrt{2 \log n}} \stackrel{\textrm{weak}}{\rightarrow} 0
\end{equation}
From the fourth moment given above we know that $\P(d_i<\frac{i}{3})=O(\frac{1}{i^2})$.  Since $s$ diverges to infinity, it follows immediately that $\sum_{i=s}^n \epsilon_i^{''}$ is almost surely zero.  Combining this with \eqref{eqn:gausstrunc}, we see that the $\epsilon^{''}$ portion of our truncation error is also negligible.

The theorem follows from our bounds on each term in the Taylor expansion.
\end{document}